\documentclass[12pt,intlimits,reqno]{article}

\usepackage{amssymb}
\usepackage{epsfig}
\usepackage{graphicx}
\usepackage{lscape}
\usepackage{amsmath}
\usepackage{amsfonts}
\usepackage{amscd}
\usepackage{cite}
\usepackage{amsthm}

\usepackage[all,cmtip]{xy}

\newenvironment{theo}[1]{\vskip+0.4cm {\fontsize{14}{16pt}\selectfont \textbf{Theorem #1.\;}}}{\vskip+0.4cm}
\newenvironment{cor}[1]{\vskip+0.4cm {\fontsize{14}{16pt}\selectfont \textbf{Corollary #1.\;}}}{\vskip+0.4cm}
\newenvironment{property}[1]{\vskip+0.4cm {\fontsize{14}{16pt}\selectfont \textbf{Property #1.\;}}}{\vskip+0.4cm}

\begin{document}

\fontsize{13}{14pt}\selectfont

\title{\textbf{Perfect Cuboid and Congruent Number Equation Solutions}}

\author{\textbf{Mamuka Meskhishvili}}
\date{}
\maketitle

\begin{center}
    {\fontsize{14}{16pt}\selectfont \textbf{Abstract} }
\end{center}
\vskip+0.5cm

A perfect cuboid (PC) is a rectangular parallelepiped with rational sides {\fontsize{14}{16pt}\selectfont $a$}, {\fontsize{14}{16pt}\selectfont $b$}, {\fontsize{14}{16pt}\selectfont $c$} whose face diagonals {\fontsize{14}{16pt}\selectfont $d_{ab}$}, {\fontsize{14}{16pt}\selectfont $d_{bc}$}, {\fontsize{14}{16pt}\selectfont $d_{ac}$} and space (body) diagonal {\fontsize{14}{16pt}\selectfont $d_s$} are rationals. The existence or otherwise of PC is a problem known since at least the time of Leonhard Euler. This research establishes equivalent conditions of PC by nontrivial rational solutions {\fontsize{14}{16pt}\selectfont $(X,Y)$} and {\fontsize{14}{16pt}\selectfont $(Z,W)$} of congruent number equation
\vskip+0.01cm
{\fontsize{14}{16pt}\selectfont
$$  y^2=x^3-N^2x,      $$
}
\vskip+0.01cm
\noindent
where product {\fontsize{14}{16pt}\selectfont $XZ$} is a square. By using such pair of solutions five pa\-ra\-met\-ri\-za\-ti\-ons of nearly-perfect cuboid (NPC) (only one face diagonal is irrational) and five equivalent conditions for PC were found. Each parametrization gives all possible NPC. For example, by using one of them -- invariant pa\-ra\-met\-ri\-za\-ti\-on for sides and diagonals of NPC are obtained:
\vskip+0.01cm
{\fontsize{14}{16pt}\selectfont
\begin{gather*}
    a=2XZN, \;\;\; b=|YW|, \;\;\; c=|X-Z|\sqrt{XZ}\,N, \\[0.2cm]
    d_{bc}=|XZ-N^2|\sqrt{XZ}, \;\;\; d_{ac}=|X+Z|\sqrt{XZ}\,N, \\[0.2cm]
    d_s=(XZ+N^2)\sqrt{XZ}\,;
\end{gather*}
}
\vskip+0.01cm
\noindent
and condition of the existence of PC is the rationality of
\vskip+0.01cm
{\fontsize{14}{16pt}\selectfont
$$  d_{ab}=\sqrt{Y^2W^2+4N^2X^2Z^2}\,.      $$
}

Because each parametrization is complete, inverse problem is discussed. For given NPC is found corresponding congruent number equation (i.e. congruent number) and its solutions.

\medskip


\medskip

\textbf{2010 AMS Classification.}
    11D25, 11D41, 11D72, 14G05, 14H52.

\vskip+2cm

\section{Introduction}
\vskip+0.5cm

Perfect cuboid (PC) problem is equivalent to the system of Diophantine equations:
\vskip+0.02cm
{\fontsize{14}{16pt}\selectfont
\begin{equation}\label{eqn1-1}
\begin{aligned}
        a^2+b^2 & =d_{ab}^{\,2}, \\[0.2cm]
        b^2+c^2 & =d_{bc}^{\,2}, \\[0.2cm]
        a^2+c^2 & =d_{ac}^{\,2}, \\[0.2cm]
        a^2+b^2+c^2 & =d_s^{\,2};
\end{aligned}
\end{equation}
}
\vskip+0.2cm
\noindent
integer (rational) solution existence. If we remove the integer (rational) con\-di\-ti\-on for space diagonal, then we get Euler cuboid (EC). If among the seven:
\vskip+0.02cm
{\fontsize{14}{16pt}\selectfont
    $$  a,\;b,\;c,\;d_{ab},\;d_{bc},\;d_{ac},\;d_s     $$
 }
\vskip+0.2cm
\noindent
only one face diagonal or side is irrational, then this is called nearly-perfect cuboid (NPC). Numerous EC and NPC have been found by the help of computers {\fontsize{14}{16pt}\selectfont \cite{1, 2, 3, 4, 5, 6, 7}}. As for the search for a PC, the computer programs of many researchers in many countries have been unsuccessful. Among the recent research {\fontsize{14}{16pt}\selectfont \cite{8}} with the help of computers it was proved that there exists no PC, the smallest side of which is less than
\vskip+0.02cm
{\fontsize{14}{16pt}\selectfont
    $$  2,1\cdot 10^{10}.      $$
    }
\vskip+0.2cm

Eventually, having recorded such large number, we come to a hypothesis that there exists no PC, but for now, this hypothesis has not been proved right yet.

\vskip+1cm
\section{Perfect Cuboid First Equation}
\vskip+0.5cm

Rewrite the system \eqref{eqn1-1} of PC as follows
\vskip+0.02cm
{\fontsize{14}{16pt}\selectfont
\begin{align*}
    a^2+d_{bc}^{\,2} & =d_s^{\,2}, \\[0.2cm]
    b^2+d_{ac}^{\,2} & =d_s^{\,2}, \\[0.2cm]
    c^2+d_{ab}^{\,2} & =d_s^{\,2}.
\end{align*}
}
\vskip+0.2cm
\noindent
Divide each equation by  {\fontsize{14}{16pt}\selectfont $d_s^{\,2}$}:
\vskip+0.02cm
{\fontsize{14}{16pt}\selectfont
\begin{align*}
    \Big(\frac{a}{d_s}\Big)^2+\Big(\frac{d_{bc}}{d_s}\Big)^2 & =1, \\[0.2cm]
    \Big(\frac{b}{d_s}\Big)^2+\Big(\frac{d_{ac}}{d_s}\Big)^2 & =1, \\[0.2cm]
    \Big(\frac{c}{d_s}\Big)^2+\Big(\frac{d_{ab}}{d_s}\Big)^2 & =1.
\end{align*}
}
\vskip+0.2cm
\noindent
Use rational parametrization formulae for unit circle
{\fontsize{14}{16pt}\selectfont
$$      x^2+y^2=1,       $$
}
\vskip+0.1cm
\noindent
its positive rational solutions are:
{\fontsize{14}{16pt}\selectfont
$$  x=\Big|\frac{1-t^2}{1+t^2}\Big|, \;\;\; y=\Big|\frac{2t}{1+t^2}\Big|,       $$
}
\vskip+0.1cm
\noindent
where {\fontsize{14}{16pt}\selectfont $t$} is arbitrary nontrivial rational number
{\fontsize{14}{16pt}\selectfont
$$  t\in\mathbb{Q}\setminus\{0;\pm1\}.       $$
}
\vskip+0.1cm
\noindent
Take  {\fontsize{14}{16pt}\selectfont $\alpha_1$}, {\fontsize{14}{16pt}\selectfont $\beta_1$}, {\fontsize{14}{16pt}\selectfont $\gamma_1$} for parametrization variables, then we obtain a system:
\vskip+0.05cm
{\fontsize{14}{16pt}\selectfont
\begin{equation}\label{eqn2-1}
\begin{gathered}
        \frac{a}{d_s}=\Big|\frac{2\alpha_1}{1+\alpha_1^2}\Big|, \;\;\; \frac{d_{bc}}{d_s}=\Big|\frac{1-\alpha_1^2}{1+\alpha_1^2}\Big|, \\[0.2cm]
        \frac{b}{d_s}=\Big|\frac{1-\beta_1^2}{1+\beta_1^2}\Big|, \;\;\; \frac{d_{ac}}{d_s}=\Big|\frac{2\beta_1}{1+\beta_1^2}\Big|, \\[0.2cm]
        \frac{c}{d_s}=\Big|\frac{2\gamma_1}{1+\gamma_1^2}\Big|, \;\;\; \frac{d_{ab}}{d_s}=\Big|\frac{1-\gamma_1^2}{1+\gamma_1^2}\Big|.
\end{gathered}
\end{equation}
}
\vskip+0.3cm
\noindent
Insert these expressions into the system \eqref{eqn1-1}:
\vskip+0.05cm
\allowdisplaybreaks
{\fontsize{14}{16pt}\selectfont
\begin{align*}
& \begin{aligned}
    \Big(\frac{2\alpha_1}{1+\alpha_1^2}\Big)^2+\Big(\frac{1-\beta_1^2}{1+\beta_1^2}\Big)^2 & =
            \Big(\frac{1-\gamma_1^2}{1+\gamma_1^2}\Big)^2, \\[0.2cm]
    \Big(\frac{1-\beta_1^2}{1+\beta_1^2}\Big)^2+\Big(\frac{2\gamma_1}{1+\gamma_1^2}\Big)^2 & =
            \Big(\frac{1-\alpha_1^2}{1+\alpha_1^2}\Big)^2, \\[0.2cm]
    \Big(\frac{2\alpha_1}{1+\alpha_1^2}\Big)^2+\Big(\frac{2\gamma_1}{1+\gamma_1^2}\Big)^2 & =
            \Big(\frac{2\beta_1}{1+\beta_1^2}\Big)^2,
\end{aligned} \\[0.2cm]
    & \;\Big(\frac{2\alpha_1}{1+\alpha_1^2}\Big)^2+\Big(\frac{1-\beta_1^2}{1+\beta_1^2}\Big)^2+\Big(\frac{2\gamma_1}{1+\gamma_1^2}\Big)^2=1.
\end{align*}
}
\vskip+0.3cm
\noindent
By using elementary properties, each equation of a system is reduced to third equation. So,

\begin{theo}{1}
\textit{The existence of PC is equivalent to the existence of non\-tri\-vi\-al rational solution of equation:}
\vskip+0.05cm
{\fontsize{14}{16pt}\selectfont
$$    \Big(\frac{2\alpha_1}{1+\alpha_1^2}\Big)^2+\Big(\frac{2\gamma_1}{1+\gamma_1^2}\Big)^2=\Big(\frac{2\beta_1}{1+\beta_1^2}\Big)^2.   $$
}
\end{theo}

\vskip+1cm
\section{Perfect Cuboid Second Equation}
\vskip+0.5cm

Rewrite the system \eqref{eqn1-1} of PC as follows
\vskip+0.02cm
{\fontsize{14}{16pt}\selectfont
\begin{align*}
    d_s^{\,2}-d_{bc}^2 & =a^2, \\[0.2cm]
    d_{ac}^{\,2}-c^2 & =a^2, \\[0.2cm]
    d_{ab}^{\,2}-b^2 & =a^2.
\end{align*}
}
\vskip+0.2cm
\noindent
Divide each equation by {\fontsize{14}{16pt}\selectfont $a^2$}:
\vskip+0.02cm
{\fontsize{14}{16pt}\selectfont
\allowdisplaybreaks
\begin{align*}
    \Big(\frac{d_s}{a}\Big)^2-\Big(\frac{d_{bc}}{a}\Big)^2 & =1, \\[0.2cm]
    \Big(\frac{d_{ac}}{a}\Big)^2-\Big(\frac{c}{a}\Big)^2 & =1, \\[0.2cm]
    \Big(\frac{d_{ab}}{a}\Big)^2-\Big(\frac{b}{a}\Big)^2 & =1.
\end{align*}
}
\vskip+0.2cm
\noindent
Use rational parametrization formulae for unit hyperbola
\vskip+0.02cm
{\fontsize{14}{16pt}\selectfont
$$  x^2-y^2=1,       $$
}
\vskip+0.2cm
\noindent
its positive rational solutions are:
\vskip+0.05cm
{\fontsize{14}{16pt}\selectfont
$$  x=\Big|\frac{1+t^2}{1-t^2}\Big|, \;\;\; y=\Big|\frac{2t}{1-t^2}\Big|,        $$
}
\vskip+0.2cm
\noindent
where {\fontsize{14}{16pt}\selectfont $t$} is arbitrary nontrivial rational number.

\noindent Take {\fontsize{14}{16pt}\selectfont $\alpha_2$}, {\fontsize{14}{16pt}\selectfont $\beta_2$} and {\fontsize{14}{16pt}\selectfont $\gamma_2$} for parametrization variables, then we obtain the system:
\vskip+0.05cm
{\fontsize{14}{16pt}\selectfont
\begin{equation}\label{eqn3-1}
\begin{gathered}
    \frac{d_s}{a}=\Big|\frac{1+\alpha_2^2}{1-\alpha_2^2}\Big|, \;\;\; \frac{d_{bc}}{a}=\Big|\frac{2\alpha_2}{1-\alpha_2^2}\Big|, \\[0.2cm]
    \frac{d_{ac}}{a}=\Big|\frac{1+\beta_2^2}{1-\beta_2^2}\Big|, \;\;\; \frac{c}{a}=\Big|\frac{2\beta_2}{1-\beta_2^2}\Big|, \\[0.2cm]
    \frac{d_{ab}}{a}=\Big|\frac{1+\gamma_2^2}{1-\gamma_2^2}\Big|, \;\;\; \frac{b}{a}=\Big|\frac{2\gamma_2}{1-\gamma_2^2}\Big|.
\end{gathered}
\end{equation}
}
\vskip+0.3cm
\noindent
Insert these expressions into system \eqref{eqn1-1}:
\vskip+0.05cm
{\fontsize{14}{16pt}\selectfont
\begin{align*}
    1+\Big(\frac{2\gamma_2}{1-\gamma_2^2}\Big)^2 & =\Big(\frac{1+\gamma_2^2}{1-\gamma_2^2}\Big)^2, \\[0.2cm]
    \Big(\frac{2\gamma_2}{1-\gamma_2^2}\Big)^2+\Big(\frac{2\beta_2}{1-\beta_2^2}\Big)^2 & =
                \Big(\frac{2\alpha_2}{1-\alpha_2^2}\Big)^2, \\[0.2cm]
    1+\Big(\frac{2\beta_2}{1-\beta_2^2}\Big)^2 & =\Big(\frac{1+\beta_2^2}{1-\beta_2^2}\Big)^2, \\[0.2cm]
    1+\Big(\frac{2\gamma_2}{1-\gamma_2^2}\Big)^2+\Big(\frac{2\beta_2}{1-\beta_2^2}\Big)^2 & =
                \Big(\frac{1+\alpha_2^2}{1-\alpha_2^2}\Big)^2.
\end{align*}
}
\vskip+0.3cm
\noindent
The first and the third equations are identities, whereas the second and fourth are equivalent. So,

\begin{theo}{2}
\textit{The existence of PC is equivalent to the existence of non\-tri\-vi\-al rational solution of equation:}
\vskip+0.05cm
{\fontsize{14}{16pt}\selectfont
$$    \Big(\frac{2\gamma_2}{1-\gamma_2^2}\Big)^2+\Big(\frac{2\beta_2}{1-\beta_2^2}\Big)^2=
            \Big(\frac{2\alpha_2}{1-\alpha_2^2}\Big)^2.      $$
}
\end{theo}

\vskip+1cm
\section{Perfect Cuboid Third Equation}
\vskip+0.5cm

Discuss second type rational parametrization for unit hyperbola
\vskip+0.05cm
{\fontsize{14}{16pt}\selectfont
$$  x^2-y^2=1,       $$
}
\vskip+0.3cm
\noindent
its positive rational solutions are:
\vskip+0.05cm
{\fontsize{14}{16pt}\selectfont
$$  x=\Big|\frac{1+t^2}{2t}\Big|, \;\;\; y=\Big|\frac{1-t^2}{2t}\Big|,        $$
}
\vskip+0.3cm
\noindent
where {\fontsize{14}{16pt}\selectfont $t$} is arbitrary nontrivial rational number.

\noindent Take {\fontsize{14}{16pt}\selectfont $\alpha_3$}, {\fontsize{14}{16pt}\selectfont $\beta_3$} and {\fontsize{14}{16pt}\selectfont $\gamma_3$}, for parametrization variables, then we obtain the system:
\vskip+0.05cm
{\fontsize{14}{16pt}\selectfont
\begin{equation}\label{eqn4-1}
\begin{gathered}
    \frac{d_s}{a}=\Big|\frac{1+\alpha_3^2}{2\alpha_3}\Big|, \;\;\; \frac{d_{bc}}{a}=\Big|\frac{1-\alpha_3^2}{2\alpha_3}\Big|, \\[0.2cm]
    \frac{d_{ac}}{a}=\Big|\frac{1+\beta_3^2}{2\beta_3}\Big|, \;\;\; \frac{c}{a}=\Big|\frac{1-\beta_3^2}{2\beta_3}\Big|, \\[0.2cm]
    \frac{d_{ab}}{a}=\Big|\frac{1+\gamma_3^2}{2\gamma_3}\Big|, \;\;\; \frac{b}{a}=\Big|\frac{1-\gamma_3^2}{2\gamma_3}\Big|.
\end{gathered}
\end{equation}
}
\vskip+0.3cm
\noindent
Insert these expressions into the system \eqref{eqn1-1}:
\vskip+0.05cm
{\fontsize{14}{16pt}\selectfont
\begin{align*}
    1+\Big(\frac{1-\gamma_3^2}{2\gamma_3}\Big)^2 & =\Big(\frac{1+\gamma_3^2}{2\gamma_3}\Big)^2, \\[0.2cm]
    \Big(\frac{1-\gamma_3^2}{2\gamma_3}\Big)^2+\Big(\frac{1-\beta_3^2}{2\beta_3}\Big)^2 & =
            \Big(\frac{1-\alpha_3^2}{2\alpha_3}\Big)^2, \\[0.2cm]
    1+\Big(\frac{1-\beta_3^2}{2\beta_3}\Big)^2 & =\Big(\frac{1+\beta_3^2}{2\beta_3}\Big)^2, \\[0.2cm]
    1+\Big(\frac{1-\gamma_3^2}{2\gamma_3}\Big)^2+\Big(\frac{1-\beta_3^2}{2\beta_3}\Big)^2 & =\Big(\frac{1+\alpha_3^2}{2\alpha_3}\Big)^2.
\end{align*}
}
\vskip+0.3cm
\noindent
The first and the third equations are identities, whereas the second and fourth are equivalent. So,

\begin{theo}{3}
\textit{The existence of PC is equivalent to the existence of non\-tri\-vi\-al rational solution of the equation:}
\vskip+0.02cm
{\fontsize{14}{16pt}\selectfont
$$    \Big(\frac{1-\gamma_3^2}{2\gamma_3}\Big)^2+\Big(\frac{1-\beta_3^2}{2\beta_3}\Big)^2=\Big(\frac{1-\alpha_3^2}{2\alpha_3}\Big)^2.    $$
}
\vskip+0.2cm
\end{theo}

Theorem 3 and Theorem 2 equations are birrationally equivalent over
\vskip+0.05cm
{\fontsize{14}{16pt}\selectfont
$$  \mathbb{Q}\setminus\{0;\pm1\}.       $$
}
\vskip+0.3cm
\noindent
Equivalency is given by rational transformation
\vskip+0.05cm
{\fontsize{14}{16pt}\selectfont
$$  \gamma_2=\frac{1-\gamma_3}{1+\gamma_3}\,, \;\;\; \beta_2=\frac{1-\beta_3}{1+\beta_3}\,, \;\;\;
            \alpha_2=\frac{1-\alpha_3}{1+\alpha_3}\,.       $$
}

\vskip+1cm
\section{Congruent Number Equation Solutions \\ Properties}
\vskip+0.5cm

By using solutions of congruent number equation (congruent curve)
\vskip+0.05cm
{\fontsize{14}{16pt}\selectfont
$$  C_N:\;\;y^2=x^3-N^2x     $$
}
\vskip+0.2cm
\noindent
it is possible to construct rational right triangles with area {\fontsize{14}{16pt}\selectfont $N$} and 3-term arithmetical progressions of squares with common difference {\fontsize{14}{16pt}\selectfont $N$} {\fontsize{14}{16pt}\selectfont \cite{9, 10}}. In both cases there is one-to-one correspondence between above mentioned sets and only points of congruent curve which are obtained by drawn tangent line through some points.

The usage of solutions of congruent number equation is not limited to the mentioned two cases. This research found the third usage.

By two solutions of congruent number equation every NPC is constructed and PC existence equivalency condition is found. It is impossible to choose arbitrary pair of solutions, they must satisfy the following condition -- the product of solutions must be a square.

First of all prove that arbitrary congruent number equation has such in\-fi\-ni\-te\-ly many solutions.

Denote the addition operation of rational points of {\fontsize{14}{16pt}\selectfont $C_N$} elliptic curve (con\-gru\-ent curve is a case of elliptic curve) by {\fontsize{14}{16pt}\selectfont $\oplus$} symbol, then the point which are obtained by drawn tangent line through {\fontsize{14}{16pt}\selectfont $P(X,Y)$} is
\vskip+0.01cm
{\fontsize{14}{16pt}\selectfont
$$  P\oplus P=2P,       $$
}
\vskip+0.2cm
\noindent
then the first coordinate of {\fontsize{14}{16pt}\selectfont $2P$} is
\vskip+0.05cm
{\fontsize{14}{16pt}\selectfont
$$  [2P]_x=\Big(\frac{X^2+N^2}{2Y}\Big)^2.      $$
}
\vskip+0.3cm
\noindent
Meaning that by drawing tangent line all first coordinates for followings
\vskip+0.02cm
{\fontsize{14}{16pt}\selectfont
$$  2P\oplus 2P=4P,\;\;3P\oplus 3P=6P,\;\;4P\oplus 4P=8P,\;\ldots       $$
}
\vskip+0.2cm
\noindent
are squares.

By drawing a secant line through {\fontsize{14}{16pt}\selectfont $P$} and {\fontsize{14}{16pt}\selectfont $2P$} points we obtain
\vskip+0.05cm
{\fontsize{14}{16pt}\selectfont
$$  P\oplus 2P=3P,       $$
}
\vskip+0.3cm
\noindent
the first coordinates of which satisfy:
\vskip+0.05cm
{\fontsize{14}{16pt}\selectfont
$$  [P]_x\cdot[2P]_x\cdot[3P]_x=d^{\,2},        $$
}
\vskip+0.3cm
\noindent
where {\fontsize{14}{16pt}\selectfont $d$} is a {\fontsize{14}{16pt}\selectfont $y$}-intercept of secant line
\vskip+0.05cm
{\fontsize{14}{16pt}\selectfont
$$  d=\frac{X\cdot[2P]_y-Y\cdot[2P]_x}{X-[2P]_x}\,.       $$
}
\vskip+0.3cm
\noindent
Meaning that the product \;{\fontsize{14}{16pt}\selectfont $[P]_x\cdot[3P]_x$}\; is a square.

By the consequent use of the reasoning we obtain

\begin{property}{1}
\textit{For arbitrary point {\fontsize{14}{16pt}\selectfont $P$} of congruent curve {\fontsize{14}{16pt}\selectfont $C_N$}, the pro\-duct
\vskip+0.05cm
{\fontsize{14}{16pt}\selectfont
$$  [kP]_x\cdot[mP]_x       $$
}
\vskip+0.3cm
\noindent
is a square, if {\fontsize{14}{16pt}\selectfont $k$} and {\fontsize{14}{16pt}\selectfont $m$} numbers have the same parity.}
\end{property}

So, congruent number equation has infinitely many rational pair solutions, the product of which is a square.

By using the solutions {\fontsize{14}{16pt}\selectfont $X$} and {\fontsize{14}{16pt}\selectfont $Z$} of congruent curves we can find all nontrivial rational solutions for special type of Kummer's surface {\fontsize{14}{16pt}\selectfont \cite{11}}.

\begin{property}{2}
\textit{All nontrivial rational solutions of the equation
\vskip+0.05cm
{\fontsize{14}{16pt}\selectfont
$$  \eta^2=\xi\zeta(\xi^2-1)(\zeta^2-1),     $$
}
\vskip+0.3cm
\noindent
are obtained by formulae
\vskip+0.05cm
{\fontsize{14}{16pt}\selectfont
$$  (\xi,\zeta,\eta)=\Big(\frac{X}{N}\,,\frac{Z}{N}\,,\frac{YW}{N^3}\Big),       $$
}
\vskip+0.3cm
\noindent
where {\fontsize{14}{16pt}\selectfont $(X,Y)$} and {\fontsize{14}{16pt}\selectfont $(Z,W)$} are arbitrary nontrivial different rational so\-lu\-ti\-ons of arbitrary {\fontsize{14}{16pt}\selectfont $C_N$} congruent number equation.}
\end{property}

Nontrivial solutions of congruent equation are the solutions with {\fontsize{14}{16pt}\selectfont $y\neq 0$} and Kummer's surface equation nontrivial solutions are the solutions with {\fontsize{14}{16pt}\selectfont $\eta\neq 0$}, {\fontsize{14}{16pt}\selectfont $\xi\neq \zeta$}.

The points obtained by drawing a secant line through point {\fontsize{14}{16pt}\selectfont $P(X,Y)$} of congruent curve {\fontsize{14}{16pt}\selectfont $C_N$} and trivial {\fontsize{14}{16pt}\selectfont $(0,0)$} and {\fontsize{14}{16pt}\selectfont $(N,0)$} points are called first and second reflected points
\vskip+0.05cm
{\fontsize{14}{16pt}\selectfont
\begin{align*}
    (X,Y) & \longrightarrow \Big(-\frac{N^2}{X}\,,-\frac{N^2Y}{X^2}\Big), \\[0.2cm]
    (X,Y) & \longrightarrow \Big(\frac{N(X+N)}{X-N}\,,\frac{2N^2Y}{(X-N)^2}\Big).
\end{align*}
}
\vskip+0.3cm
\noindent
By drawing a secant line through the third trivial point {\fontsize{14}{16pt}\selectfont $(-N,0)$}, we obtain point
\vskip+0.05cm
{\fontsize{14}{16pt}\selectfont
$$  \Big(\frac{N(N-X)}{X+N}\,,\frac{2N^2Y}{(X+N)^2}\Big),     $$
}
\vskip+0.3cm
\noindent
which is the result of the composition of the first and second reflected trans\-for\-ma\-ti\-ons.

\begin{property}{3}
If {\fontsize{14}{16pt}\selectfont $X$} and {\fontsize{14}{16pt}\selectfont $Z$} are two solutions of congruent number equation the product of \textit{which is a square, then the product of the first reflected solutions and second reflected solutions
\vskip+0.05cm
{\fontsize{14}{16pt}\selectfont
$$  \Big(-\frac{N^2}{X}\Big)\cdot\Big(-\frac{N^2}{Z}\Big) \;\;\text{\fontsize{12}{16pt}\selectfont and}\;\;
            \frac{N(X+N)}{X-N}\cdot \frac{N(Z+N)}{Z-N}      $$
}
\vskip+0.3cm
\noindent
are also squares.}
\end{property}

Indeed,
\vskip+0.05cm
{\fontsize{14}{16pt}\selectfont
$$  \frac{N(X+N)}{X-N}\cdot \frac{N(Z+N)}{Z-N}=\frac{N^2Y^2W^2}{(X-N)^2(Z-N)^2XZ}\,.        $$
}

\vskip+1cm
\section{First and its Reflected Parametrizations of NPC. PC Conditions}
\vskip+0.5cm

From the Theorem 1
\vskip+0.05cm
{\fontsize{14}{16pt}\selectfont
\begin{equation}\label{eqn6-1}
    \Big(\frac{\gamma_1}{1+\gamma_1^2}\Big)^2=\frac{\alpha_1^2\Big(1-\Big(\dfrac{\alpha_1}{\beta_1}\Big)^2\Big)(1-(\alpha_1\beta_1)^2)}
                {\Big(\dfrac{\alpha_1}{\beta_1}\Big)^2(1+\alpha_1^2)^2(1+\beta_1^2)^2}\,.
\end{equation}
}
\vskip+0.3cm
\noindent
Denoted the numerator by {\fontsize{14}{16pt}\selectfont $\eta_1^2$}
\vskip+0.05cm
{\fontsize{14}{16pt}\selectfont
$$  \eta_1^2=\Big(\frac{\alpha_1}{\beta_1}\Big)(\alpha_1\beta_1)\Big(1-\Big(\frac{\alpha_1}{\beta_1}\Big)^2\Big)
           \big(1-(\alpha_1\beta_1)^2\big),    $$
}
\vskip+0.3cm
\noindent
and using Property 2:
\vskip+0.05cm
{\fontsize{14}{16pt}\selectfont
$$  \alpha_1\beta_1=\frac{X_1}{N_1}\,, \;\;\; \frac{\alpha_1}{\beta_1}=\frac{Z_1}{N_1}\,. $$
}
\vskip+0.3cm
\noindent
So
\vskip+0.05cm
{\fontsize{14}{16pt}\selectfont
$$  \alpha_1=\frac{\sqrt{X_1Z_1}}{N_1}\,, \;\;\; \beta_1=\sqrt{\frac{X_1}{Z_1}}\,, \;\;\; \eta_1=\frac{Y_1W_1}{N_1^3}\,,   $$
}
\vskip+0.3cm
\noindent
where {\fontsize{14}{16pt}\selectfont $(X_1,Y_1)$} and {\fontsize{14}{16pt}\selectfont $(Z_1,W_1)$} are arbitrary nontrivial different solutions of arbitrary {\fontsize{14}{16pt}\selectfont $C_{N_1}$} congruent number equation. By using Property 1 the solutions {\fontsize{14}{16pt}\selectfont $X_1$} and {\fontsize{14}{16pt}\selectfont $Z_1$}, the product (ratio) of which is a square, are infinitely many. For these solutions the equation \eqref{eqn6-1} is:
\vskip+0.05cm
{\fontsize{14}{16pt}\selectfont
$$  \frac{\gamma_1}{1+\gamma_1^2}=\frac{Y_1W_1}{(X_1Z_1+N_1^2)(X_1+Z_1)}\,.       $$
}
\vskip+0.3cm

If in the last expression {\fontsize{14}{16pt}\selectfont $\gamma_1$} is rational, we obtain PC.

Insert the given expressions for variables {\fontsize{14}{16pt}\selectfont $\alpha_1$}, {\fontsize{14}{16pt}\selectfont $\beta_1$} and {\fontsize{14}{16pt}\selectfont $\gamma_1$} in system \eqref{eqn2-1}. We obtain the first parametrization formulae for sides and diagonals.

\begin{theo}{4}
\textit{Complete parametrization of NPC is given by formulae:
\vskip+0.05cm
{\fontsize{14}{16pt}\selectfont
\begin{align*}
    a & =\frac{2N_1\sqrt{X_1Z_1}}{X_1Z_1+N_1^2}\cdot d_s, \\[0.2cm]
    b & =\Big|\frac{X_1-Z_1}{X_1+Z_1}\Big|\cdot d_s, \\[0.2cm]
    c & =\Big|\frac{2Y_1W_1}{(X_1Z_1+N_1^2)(X_1+Z_1)}\Big|\cdot d_s, \\[0.2cm]
    d_{ac} & =\frac{2\sqrt{X_1Z_1}}{|X_1+Z_1|}\cdot d_s, \\[0.2cm]
    d_{bc} & =\frac{|X_1Z_1-N_1^2|}{X_1Z_1+N_1^2}\cdot d_s,
\end{align*}
}
\vskip+0.3cm
\noindent
where {\fontsize{14}{16pt}\selectfont $(X_1,Y_1)$} and {\fontsize{14}{16pt}\selectfont $(Z_1,W_1)$} are arbitrary nontrivial different rational so\-lu\-ti\-ons of arbitrary {\fontsize{14}{16pt}\selectfont $C_{N_1}$} congruent number equation, the product {\fontsize{14}{16pt}\selectfont $X_1Z_1$} is a square. }

\textit{Con\-di\-ti\-on of the existence of PC is the rationality of}
\vskip+0.05cm
{\fontsize{14}{16pt}\selectfont
$$  d_{ab}=\sqrt{1-\Big(\frac{2Y_1W_1}{(X_1Z_1+N_1^2)(X_1+Z_1)}\Big)^2}\cdot d_s.     $$
}
\end{theo}

Based on Property 3 solutions {\fontsize{14}{16pt}\selectfont $(X_1,Y_1)$} and {\fontsize{14}{16pt}\selectfont $(Z_1,W_1)$} can be replaced by the first and second reflected points. Obtained formulae do not change by the first reflected transformation, while they are changed by second reflected transformation. As a result we obtain reflected parametrization.

\begin{cor}{1}
\textit{Complete parametrization of NPC is given by formulae:
\vskip+0.05cm
{\fontsize{14}{16pt}\selectfont
\begin{align*}
    a' & =\frac{|Y_1W_1|}{(X_1Z_1+N_1^2)\sqrt{X_1Z_1}}\cdot d_s{}{\!'}, \\[0.2cm]
    b' & =\Big|\frac{N_1(Z_1-X_1)}{X_1Z_1-N_1^2}\Big|\cdot d_s{}{\!'}, \\[0.2cm]
    c' & =\Big|\frac{2N_1Y_1W_1}{X_1^2Z_1^2-N_1^4}\Big|\cdot d_s{}{\!'}, \\[0.2cm]
    d_{ac}{}{\!\!\!'} & =\Big|\frac{Y_1W_1}{(X_1Z_1-N_1^2)\sqrt{X_1Z_1}}\Big|\cdot d_s{}{\!'}, \\[0.2cm]
    d_{bc}{}{\!\!\!'} & =\frac{|N_1(X_1+Z_1)|}{X_1Z_1+N_1^2}\cdot d_s{}{\!'},
\end{align*}
}
\vskip+0.3cm
\noindent
where {\fontsize{14}{16pt}\selectfont $(X_1,Y_1)$} and {\fontsize{14}{16pt}\selectfont $(Z_1,W_1)$} are arbitrary nontrivial different rational so\-lu\-ti\-ons of arbitrary {\fontsize{14}{16pt}\selectfont $C_{N_1}$} congruent number equation, the product {\fontsize{14}{16pt}\selectfont $X_1Z_1$} is a square. }

\textit{Con\-di\-ti\-on of the existence of PC is the rationality of}
\vskip+0.05cm
{\fontsize{14}{16pt}\selectfont
$$  d_{ab}{}{\!\!\!'}=\sqrt{1-\Big(\frac{2N_1Y_1W_1}{X_1^2Z_1^2-N_1^4}\Big)^2}\cdot d_s{}{\!'}.    $$
}
\end{cor}

\vskip+1cm
\section{Second and its Reflected Parametrizations of NPC. PC Conditions}
\vskip+0.5cm

From the Theorem 2
\vskip+0.05cm
{\fontsize{14}{16pt}\selectfont
\begin{equation}\label{eqn7-1}
    \Big(\frac{\gamma_2}{1-\gamma_2^2}\Big)^2=\frac{\beta_2^2\Big(1-\Big(\dfrac{\beta_2}{\alpha_2}\Big)^2\Big)(1-(\alpha_2\beta_2)^2)}
                {\dfrac{\beta_2^2}{\alpha_2^2}\,(1-\alpha_2^2)^2(1-\beta_2^2)^2}\,.
\end{equation}
}
\vskip+0.3cm
\noindent
Denote the numerator by {\fontsize{14}{16pt}\selectfont $\eta_2^2$}
\vskip+0.05cm
{\fontsize{14}{16pt}\selectfont
$$  \eta_2^2=\Big(\frac{\beta_2}{\alpha_2}\Big)(\beta_2\alpha_2)\Big(1-\Big(\frac{\beta_2}{\alpha_2}\Big)^2\Big)
                \big(1-(\alpha_2\beta_2)^2\big),   $$
}
\vskip+0.3cm
\noindent
and using Property 2
\vskip+0.05cm
{\fontsize{14}{16pt}\selectfont
$$  \frac{\beta_2}{\alpha_2}=\frac{X_2}{N_2}\,, \;\;\; \beta_2\alpha_2=\frac{Z_2}{N_2}      $$
}
\vskip+0.3cm
\noindent
where {\fontsize{14}{16pt}\selectfont $(X_2,Y_2)$} and {\fontsize{14}{16pt}\selectfont $(Z_2,W_2)$} are arbitrary nontrivial different solutions of arbitrary {\fontsize{14}{16pt}\selectfont $C_{N_2}$} congruent number equation. So,
\vskip+0.05cm
{\fontsize{14}{16pt}\selectfont
$$  \alpha_2=\sqrt{\frac{Z_2}{X_2}}\,, \;\;\; \beta_2=\frac{\sqrt{X_2Z_2}}{N_2}\,, \;\;\; \eta_2=\frac{Y_2W_2}{N_2^3}\,.        $$
}
\vskip+0.3cm
\noindent
Based on Property 1, solutions the product (ratio) of which is a square are infinitely many. For such solutions the equation
\eqref{eqn7-1} is
\vskip+0.05cm
{\fontsize{14}{16pt}\selectfont
$$  \frac{\gamma_2}{1-\gamma_2^2}=\frac{Y_2W_2}{(X_2-Z_2)(N_2^2-X_2Z_2)}\,.       $$
}
\vskip+0.3cm
\noindent
If {\fontsize{14}{16pt}\selectfont $\gamma_2$} is rational, then we obtain PC.

Insert the obtained parametrization variables {\fontsize{14}{16pt}\selectfont $\alpha_2$}, {\fontsize{14}{16pt}\selectfont $\beta_2$} and {\fontsize{14}{16pt}\selectfont $\gamma_2$} into system \eqref{eqn3-1} we obtain third parametrization formulae for sides and diagonals.

\begin{theo}{5}
\textit{Complete parametrization of NPC is given by formulae:
\vskip+0.05cm
{\fontsize{14}{16pt}\selectfont
\begin{align*}
    b & =\Big|\frac{2Y_2W_2}{(X_2-Z_2)(N_2^2-X_2Z_2)}\Big|\cdot a, \\[0.2cm]
    c & =\frac{2N_2\sqrt{X_2Z_2}}{|N_2^2-X_2Z_2|}\cdot a, \\[0.2cm]
    d_{bc} & =\frac{2\sqrt{X_2Z_2}}{|X_2-Z_2|}\cdot a, \\[0.2cm]
    d_{ac} & =\frac{N_2^2+X_2Z_2}{|N_2^2-X_2Z_2|}\cdot a, \\[0.2cm]
    d_s & =\Big|\frac{X_2+Z_2}{X_2-Z_2}\Big|\cdot a,
\end{align*}
}
\vskip+0.3cm
\noindent
where {\fontsize{14}{16pt}\selectfont $(X_2,Y_2)$} and {\fontsize{14}{16pt}\selectfont $(Z_2,W_2)$} are arbitrary nontrivial different rational so\-lu\-ti\-ons of arbitrary {\fontsize{14}{16pt}\selectfont $C_{N_2}$} congruent number equation, the product {\fontsize{14}{16pt}\selectfont $X_2Z_2$} is a square. }

\textit{Condition of the existence of PC is the rationality of}
\vskip+0.05cm
{\fontsize{14}{16pt}\selectfont
$$  d_{ab}=\sqrt{1+\Big(\frac{2Y_2W_2}{(X_2-Z_2)(N_2^2-X_2Z_2)}\Big)^2}\cdot a.     $$
}
\end{theo}

Obtained parametrization is invariable under the first ref\-lec\-ted tran\-sfor\-ma\-ti\-on but the second ref\-lec\-ted tran\-sfor\-ma\-ti\-on gives new parametrization.

\begin{cor}{2}
\textit{Complete parametrization of NPC is given by formulae:
\vskip+0.05cm
{\fontsize{14}{16pt}\selectfont
\begin{align*}
    b' & =\Big|\frac{2Y_2W_2}{N_2(Z_2^2-X_2^2)}\Big|\cdot a', \\[0.2cm]
    c' & =\Big|\frac{Y_2W_2}{N_2(X_2+Z_2)\sqrt{X_2Z_2}}\Big|\cdot a', \\[0.2cm]
    d_{bc}{}{\!\!\!'} & =\Big|\frac{Y_2W_2}{N_2(Z_2-X_2)\sqrt{X_2Z_2}}\Big|\cdot a', \\[0.2cm]
    d_{ac}{}{\!\!\!'} & =\frac{X_2Z_2+N_2^2}{|N_2(X_2+Z_2)|}\cdot a', \\[0.2cm]
    d_s{}{\!'} & =\Big|\frac{X_2Z_2-N_2^2}{N_2(Z_2-X_2)}\Big|\cdot a',
\end{align*}
}
\vskip+0.3cm
\noindent
where {\fontsize{14}{16pt}\selectfont $(X_2,Y_2)$} and {\fontsize{14}{16pt}\selectfont $(Z_2,W_2)$} are arbitrary nontrivial different rational so\-lu\-ti\-ons of arbitrary {\fontsize{14}{16pt}\selectfont $C_{N_2}$} congruent number equation, the product {\fontsize{14}{16pt}\selectfont $X_2Z_2$} is a square. }

\textit{Condition of the existence of PC is the rationality of}
\vskip+0.05cm
{\fontsize{14}{16pt}\selectfont
$$  d_{ab}{}{\!\!\!'}=\sqrt{1+\Big(\frac{2Y_2W_2}{N_2(Z_2^2-X_2^2)}\Big)^2}\cdot a'.      $$
}
\end{cor}

\vskip+1cm
\section{Invariant Parametrization of NPC. PC Invariant Condition}
\vskip+0.5cm

From the Theorem 3
\vskip+0.05cm
{\fontsize{14}{16pt}\selectfont
\begin{equation}\label{eqn8-1}
    \Big(\frac{1-\gamma_3^2}{\gamma_3}\Big)^2=\frac{\Big(\dfrac{\alpha_3}{\beta_3}\Big)(\alpha_3\beta_3)
            \Big(1-\Big(\dfrac{\alpha_3}{\beta_3}\Big)^2\Big)(1-(\alpha_3\beta_3)^2)}{\alpha_3^4}\,.
\end{equation}
}
\vskip+0.3cm
\noindent
Denote the numerator by {\fontsize{14}{16pt}\selectfont $\eta_3^2$}
\vskip+0.05cm
{\fontsize{14}{16pt}\selectfont
$$  \eta_3^2=\Big(\frac{\alpha_3}{\beta_3}\Big)(\alpha_3\beta_3)\Big(1-\Big(\frac{\alpha_3}{\beta_3}\Big)^2\Big)
                \big(1-(\alpha_3\beta_3)^2\big).    $$
}
\vskip+0.3cm
\noindent
By Property 2
\vskip+0.05cm
{\fontsize{14}{16pt}\selectfont
$$   \frac{\alpha_3}{\beta_3}=\frac{X}{N}, \;\;\; \alpha_3\beta_3=\frac{Z}{N}\,,    $$
}
\vskip+0.3cm
\noindent
where {\fontsize{14}{16pt}\selectfont $(X,Y)$} and {\fontsize{14}{16pt}\selectfont $(Z,W)$} are arbitrary nontrivial different solutions of ar\-bit\-ra\-ry {\fontsize{14}{16pt}\selectfont $C_N$} congruent number equation. Because solutions the product (ratio) of which is a square are infinitely many (Property 1), for these solutions
\vskip+0.05cm
{\fontsize{14}{16pt}\selectfont
$$  \alpha_3=\frac{\sqrt{XZ}}{N}\,, \;\;\; \beta_3=\sqrt{\frac{Z}{X}}\,, \;\;\; \eta_3=\frac{YW}{N^3}\,,        $$
}
\vskip+0.3cm
\noindent
also from \eqref{eqn8-1}
\vskip+0.05cm
{\fontsize{14}{16pt}\selectfont
$$  \frac{1-\gamma_3^2}{\gamma_3}=\frac{YW}{XZN}\,.       $$
}
\vskip+0.3cm
\noindent
System \eqref{eqn4-1} gives the fifth parametrization formulae of sides and diagonals of NPC.

\begin{theo}{6}
\textit{Complete parametrization of NPC is given by formulae:
\vskip+0.05cm
{\fontsize{14}{16pt}\selectfont
\begin{align*}
    b & =\frac{|YW|}{2XZN}\cdot a, \\[0.2cm]
    c & =\frac{|X-Z|}{2\sqrt{XZ}}\cdot a, \\[0.2cm]
    d_{bc} & =\frac{|N^2-XZ|}{2N\sqrt{XZ}}\cdot a, \\[0.2cm]
    d_{ac} & =\frac{|X+Z|}{2\sqrt{XZ}}\cdot a, \\[0.2cm]
    d_s & =\frac{N^2+XZ}{2N\sqrt{XZ}}\cdot a,
\end{align*}
}
\vskip+0.3cm
\noindent
where {\fontsize{14}{16pt}\selectfont $(X,Y)$} and {\fontsize{14}{16pt}\selectfont $(Z,W)$} are arbitrary nontrivial different rational so\-lu\-ti\-ons of arbitrary {\fontsize{14}{16pt}\selectfont $C_N$} congruent number equation, the product {\fontsize{14}{16pt}\selectfont $XZ$} is a square. }

\textit{Condition of the existence of PC is the rationality of }
\vskip+0.05cm
{\fontsize{14}{16pt}\selectfont
$$  d_{ab}=\sqrt{1+\Big(\frac{YW}{2NXZ}\Big)^2}\cdot a.         $$
}
\end{theo}
\vskip+0.4cm

As preceeding parametrizations the last one is invariable under the first reflected transformation. Though second reflected transformation gives
\vskip+0.05cm
{\fontsize{14}{16pt}\selectfont
\begin{align*}
    a' & =\frac{|YW|}{2XZN}\cdot b', \\[0.2cm]
    c' & =\frac{|X-Z|}{2\sqrt{XZ}}\cdot b'.
\end{align*}
}
\vskip+0.3cm
\noindent
So, the second reflected transformation interchanges sides {\fontsize{14}{16pt}\selectfont $a$} and {\fontsize{14}{16pt}\selectfont $b$}, though side {\fontsize{14}{16pt}\selectfont $c$} is left unchanged. Due to this, we call the last parametrization an invariant parametrization. Invariant parametrization is the most convinient for construction of concrete NPC.

Take {\fontsize{14}{16pt}\selectfont $N=5$} and following solutions {\fontsize{14}{16pt}\selectfont \cite{12}}:
\vskip+0.05cm
{\fontsize{14}{16pt}\selectfont
\begin{gather*}
    X=\frac{25}{2^2}\,, \;\;\; Y=\frac{75}{2^3}\,, \\[0.2cm]
    Z=\frac{1681}{12^2}\,, \;\;\; W=\frac{62279}{12^3}\,,
\end{gather*}
}
\vskip+0.3cm
\noindent
By the first parametrization and its reflection:
\vskip+0.05cm
{\fontsize{14}{16pt}\selectfont
\begin{gather*}
    a=5079408, \;\;\; b=1762717, \;\;\; c=2242044, \\[0.2cm]
    d_{ac}=5552220, \;\;\; d_{bc}=2852005, \\[0.2cm]
    d_s=5825317; \\[0.5cm]
    a'=5035485, \;\;\; b'=7050868, \;\;\; c'=8968176, \\[0.2cm]
    d_{ac}{}{\!\!\!'}=10285149, \;\;\; d_{bc}{}{\!\!\!'}=11408020, \\[0.2cm]
    d_s{}{\!'}=12469925.
\end{gather*}
}
\vskip+0.3cm

By the second parametrization and its reflection:
\vskip+0.05cm
{\fontsize{14}{16pt}\selectfont
\begin{gather*}
    a=863005, \;\;\; b=2242044, \;\;\; c=1537008, \\[0.2cm]
    d_{ac}=1762717, \;\;\; d_{bc}=2718300, \\[0.2cm]
    d_s=2852005; \\[0.5cm]
    a'=8063044, \;\;\; b'=11210220, \;\;\; c'=3559017, \\[0.2cm]
    d_{ac}{}{\!\!\!'}=8813585, \;\;\; d_{bc}{}{\!\!\!'}=11761617, \\[0.2cm]
    d_s{}{\!'}=14260025.
\end{gather*}
}
\vskip+0.3cm

By invariant parametrization:
\vskip+0.05cm
{\fontsize{14}{16pt}\selectfont
\begin{gather*}
    a=9840, \;\;\; b=4557, \;\;\; c=3124, \\[0.2cm]
    d_{ac}=10324, \;\;\; d_{bc}=5525, \\[0.2cm]
    d_s=11285.
\end{gather*}
}
\vskip+0.3cm

There have been attempts to connect the solutions of congruent number equation with NPC earlier, we should mention {\fontsize{14}{16pt}\selectfont \cite{8, 11, 13, 14}}. These researches show the existence of such relationship, but due to the complexity of their formulae this relationship is not expressed explicitly. That is why NPC parametrizations and PC existence equivalency conditions have not been produced. Elliptic curves associated with PC are considered in {\fontsize{14}{16pt}\selectfont \cite{15, 16}}.

\vskip+1cm
\section{Finding Congruent Number Equation and its Solutions by Given NPC}
\vskip+0.5cm

As shown in preceeding part, using invariant parametrization NPC is con\-struc\-ted by each of four pairs of solutions of congruent number equation. Place solution pairs in square brackets.
\vskip+0.05cm
{\fontsize{14}{16pt}\selectfont
\begin{align*}
\ \hskip-1.6cm {\rm I.} &\qquad\qquad [X;Z] \longrightarrow (a,b,c), \\[0.2cm]
\ \hskip-1.6cm {\rm II.} &\qquad\qquad \Big[-\frac{N^2}{X}\,;-\frac{N^2}{Z}\Big]\longrightarrow (a,b,c), \\[0.2cm]
\ \hskip-1.6cm {\rm III.} &\qquad\qquad \Big[N\,\frac{X+N}{X-N}\,;N\,\frac{Z+N}{Z-N}\Big]\longrightarrow(b,a,c), \\[0.2cm]
\ \hskip-1.6cm {\rm IV.} &\qquad\qquad \Big[N\,\frac{N-X}{N+X}\,;N\,\frac{N-Z}{N+Z}\Big]\longrightarrow(b,a,c).
\end{align*}
}
\vskip+0.3cm
\noindent
The aim of this part is the consideration of inverse problem. For given NPC, find corresponding congruent number equation (i.e. congruent number {\fontsize{14}{16pt}\selectfont $N$}) and its four solution pairs.

From invariant parametrization (Theorem 6)
\vskip+0.05cm
{\fontsize{14}{16pt}\selectfont
\begin{align*}
    \frac{c}{a} & =\frac{|X-Z|}{2\sqrt{XZ}}\,, \\[0.2cm]
    \frac{d_s}{a} & =\frac{N^2+XZ}{2N\sqrt{XZ}}\,.
\end{align*}
}
\vskip+0.3cm
\noindent
By solving the system:
\vskip+0.05cm
{\fontsize{14}{16pt}\selectfont
\begin{gather*}
    \frac{X}{Z}=\Big(\frac{d_{ac}\pm c}{a}\Big)^2, \;\;\; \frac{XZ}{N^2}=\Big(\frac{d_s\pm d_{bc}}{a}\Big)^2; \\[0.2cm]
\begin{aligned}
    \frac{X}{N} & =\pm\,\frac{(d_{ac}\pm c)(d_s\pm d_{bc})}{a^2}\,, \\[0.2cm]
    \frac{Z}{N} & =\pm\,\frac{d_s\pm d_{bc}}{d_{ac}\pm c}\,.
\end{aligned}
\end{gather*}
}
\vskip+0.3cm
\noindent
In total there are eight solution pairs, only four of them satisfy congruent number equation solution conditions:
\vskip+0.05cm
{\fontsize{14}{16pt}\selectfont
$$  -1<\frac{X}{N}<0 \;\;\text{\fontsize{12}{16pt}\selectfont or}\;\; \frac{X}{N}>1.       $$
}
\vskip+0.3cm

Two more solution pairs are obtained by interchange of {\fontsize{14}{16pt}\selectfont $X$} and {\fontsize{14}{16pt}\selectfont $Z$} so, finally there are only two pairs of solutions connected by the first reflected transformation:
\vskip+0.05cm
{\fontsize{14}{16pt}\selectfont
\begin{align*}
\ \hskip-0.7cm {\rm I.} &\qquad\;  \frac{X}{N}=\frac{(d_{ac}+c)(d_s+d_{bc})}{a^2}\,; \;\;\; \frac{Z}{N}=\frac{d_s+d_{bc}}{d_{ac}+c}\,. \\[0.2cm]
\ \hskip-0.7cm {\rm II.} &\quad  -\frac{N}{X}=-\frac{(d_{ac}-c)(d_s-d_{bc})}{a^2}\,; \;\;\; -\frac{N}{Z}=-\frac{d_s-d_{bc}}{d_{ac}-c}\,.
\end{align*}
}
\vskip+0.3cm

Third and fourth pairs are obtained by replacing {\fontsize{14}{16pt}\selectfont $a\leftrightarrow b$}. By solving the system:
\vskip+0.05cm
{\fontsize{14}{16pt}\selectfont
\begin{align*}
    \frac{c}{b} & =\frac{|X-Z|}{2\sqrt{XZ}}\,, \\[0.2cm]
    \frac{d_s}{b} & =\frac{N^2+XZ}{2N\sqrt{XZ}}\,;
\end{align*}
}
\vskip+0.3cm
\noindent
obtain:
\vskip+0.05cm
{\fontsize{14}{16pt}\selectfont
\begin{align*}
    \quad {\rm III.} &\quad \frac{X+N}{X-N}=\frac{d_s+d_{ac}}{d_{bc}+c}\,; \;\;\; \frac{Z+N}{Z-N}=\frac{(d_{bc}+c)(d_s+d_{ac})}{b^2}\,. \\[0.4cm]
    \quad {\rm IV.} &\quad \frac{N-X}{N+X}=-\frac{d_s-d_{ac}}{d_{bc}-c}\,; \;\;\; \frac{N-Z}{N+Z}=-\frac{(d_{bc}-c)(d_s-d_{ac})}{b^2}\,.
\end{align*}
}
\vskip+0.3cm

Again III and IV pairs are connected by the first reflected transformation when, {\fontsize{14}{16pt}\selectfont ${\rm I}\leftrightarrow{\rm III}$}, and {\fontsize{14}{16pt}\selectfont ${\rm II}\leftrightarrow{\rm IV}$} are connected by the second reflected tran\-sfor\-ma\-ti\-on.

To find congruent number {\fontsize{14}{16pt}\selectfont $N$}, consider squares
\vskip+0.05cm
{\fontsize{14}{16pt}\selectfont
$$  X(X^2-N^2) \;\;\text{\fontsize{12}{16pt}\selectfont and}\;\; Z(Z^2-N^2)       $$
}
\vskip+0.3cm
\noindent
and following is a square as well:
\vskip+0.05cm
{\fontsize{14}{16pt}\selectfont
\begin{gather*}
    N^3\cdot\frac{X}{N}\,\Big(\Big(\frac{X}{N}\Big)^2-1\Big) \longrightarrow N\cdot\frac{X}{N}\,\Big(\Big(\frac{X}{N}\Big)^2-1\Big), \\[0.2cm]
    N\cdot\frac{Z}{N}\,\Big(\Big(\frac{Z}{N}\Big)^2-1\Big).
\end{gather*}
}
\vskip+0.3cm
\noindent
By these conditions, congruent number {\fontsize{14}{16pt}\selectfont $N$} is obtained.

Discuss concrete numerical case for NPC:
\vskip+0.05cm
{\fontsize{14}{16pt}\selectfont
\begin{gather*}
    a=672, \;\;\; c=104, \;\;\; b=153, \\[0.2cm]
    d_{ac}=680, \;\;\; d_{bc}=185, \\[0.2cm]
    d_s=697.
\end{gather*}
}
\vskip+0.3cm
\noindent
Receive:
\vskip+0.05cm
{\fontsize{14}{16pt}\selectfont
$$  \frac{X}{N}=\frac{49}{32} \;\;\text{\fontsize{12}{16pt}\selectfont and}\;\; \frac{Z}{N}=\frac{9}{8}\,.     $$
}
\vskip+0.3cm
\noindent
{\fontsize{14}{16pt}\selectfont $N$} is found by removing squares
\vskip+0.05cm
{\fontsize{14}{16pt}\selectfont
$$  N\cdot\frac{49}{32}\,\Big(\Big(\frac{49}{32}\Big)^2-1\Big) \longrightarrow N\cdot\frac{17}{2}\,,       $$
}
\vskip+0.3cm
\noindent
because congruent number {\fontsize{14}{16pt}\selectfont $N$} is squarefree:
\vskip+0.05cm
{\fontsize{14}{16pt}\selectfont
$$  N=34.       $$
}
\vskip+0.3cm

The same {\fontsize{14}{16pt}\selectfont $N$} is obtained by second solution. The third and fourth pairs of solutions are obtained from:
\vskip+0.05cm
{\fontsize{14}{16pt}\selectfont
$$  \frac{X+N}{X-N}=\frac{81}{17} \;\;\text{\fontsize{12}{16pt}\selectfont and}\;\; \frac{Z+N}{Z-N}=17.            $$
}
\vskip+0.3cm

Summarizing the results we conclude: NPC with sides
\vskip+0.05cm
{\fontsize{14}{16pt}\selectfont
$$  672,\;\;104,\;\;153     $$
}
\vskip+0.3cm
\noindent
is obtained by each of four pairs of solutions:
\vskip+0.05cm
{\fontsize{14}{16pt}\selectfont
\begin{align*}
    \frac{833}{16}\;\; &\text{\fontsize{12}{16pt}\selectfont and}\;\; \frac{153}{4}\,, \\[0.2cm]
    -\frac{1088}{49} \;\;&\text{\fontsize{12}{16pt}\selectfont and}\;\; -\frac{272}{9}\,, \\[0.2cm]
    162 \;\; &\text{\fontsize{12}{16pt}\selectfont and}\;\; 578, \\[0.2cm]
    -\frac{578}{81} \;\; &\text{\fontsize{12}{16pt}\selectfont and}\;\; -2;
\end{align*}
}
\vskip+0.3cm
\noindent
of following congruent number equation
\vskip+0.05cm
{\fontsize{14}{16pt}\selectfont
$$  y^2=x(x^2-34^2)     $$
}
\vskip+0.3cm
\noindent
by invariant parametrization.

The first and second parametrizations give other congruent number equ\-a\-ti\-ons and corresponding pairs of solutions. For given NPC:

by first parametrization {\fontsize{14}{16pt}\selectfont $N_1=4305$}
\vskip+0.05cm
{\fontsize{14}{16pt}\selectfont
\begin{align*}
    [X_1;Z_1] & =\Big[\frac{452025}{64}\,;\frac{18081}{4}\Big], \\[0.2cm]
    \Big[-\frac{N_1^2}{X_1}\,;-\frac{N_1^2}{Z_1}\Big] & =\big[-2624;-4100\big].
\end{align*}
}
\vskip+0.3cm

By second parametrization {\fontsize{14}{16pt}\selectfont $N_2=1717170$}
\vskip+0.05cm
{\fontsize{14}{16pt}\selectfont
\begin{align*}
    [X_2;Z_2] & =\big[165191754;3016650\big], \\[0.2cm]
    \Big[-\frac{N_2^2}{X_2}\,;-\frac{N_2^2}{Z_2}\Big] & =\big[-17850;-977466\big].
\end{align*}
}
\vskip+0.3cm

Given examples show once again that invariant parametrization is the most convinient for calculation of concrete numerical cases.

\vskip+1cm

\vskip+0.5cm

\noindent Author's address:

\medskip

\noindent {Georgian-American High School, 18 Chkondideli Str., Tbilisi 0180, Georgia.}

\noindent {\small \textit{E-mail:} \texttt{director@gahs.edu.ge} }

\end{document}